\documentclass[12pt]{article}
 \usepackage[german]{babel}
\usepackage{amssymb}
\usepackage{amsmath}
\usepackage{graphicx}
\usepackage{bbm}
\usepackage[right]{eurosym}
\usepackage[arrow, matrix, curve]{xy}
\def\vs{\vspace}
\def\noi{\noindent}

\def\IN{\mathbb N}

\def\IR{\mathbb R}

\def\IQ{\mathbb Q}

\def\an{\mathrm{an}}
\def\exp{\mathrm{exp}}

\def\ma{\mathcal}
\begin{document}
\begin{center}
{\bf \Large Logarithms, constructible functions and integration on non-archimedean models of the theory of the real field with restricted analytic functions with value group of finite archimedean rank}
\end{center}

\vspace{0.5cm} \centerline{Tobias Kaiser}

\vspace{0.7cm}\noi {\footnotesize {\bf Abstract.} 
	Given a model of the theory of the real field with restricted analytic functions such that its value group has finite archimedean rank we show how one can extend the restricted logarithm to a global logarithm with values in the polynomial ring over the model with dimension the archimedean rank. The logarithms are determined by algebraic data from the  model, namely by a section of the model and by an embedding of the value group into its Hahn group.
	If the archimedean rank of the value group coincides with the rational rank the logarithms are equivalent.
	We illustrate how one can embed such a logarithm into a model of the real field with restricted analytic functions and exponentiation. This allows us to define constructible functions with good lifting properties. As an application we establish a Lebesgue measure and integration theory with values in the polynomial ring, extending and strengthening the construction in [T. Kaiser: Lebesgue measure and integration theory on non-archimedean real closed fields with archimedean value group. Proc. Lond. Math. Soc. 116 (2018), no. 2, 209-247.].}   

\section*{Introduction}

We are interested in analysis on non-standard structures. To be more precise, we want to establish important analytic concepts as integration in situations beyond the real field.\\
A promising framework is to work in a non-archimedean model of the theory $T_\an$ of the real field with restricted analytic functions $\IR_\an$ (see Van den Dries et al. [6] for deep results on this theory) since $\IR_\an$ defines locally all analytic functions but show similarities to algebraic geometry as Puiseux series expansion and polynomially boundedness (see for example [12, 13] for results on analytic functions in models of $T_\an$).\\
Yet for integration, the o-minimal structure $\IR_\an$ is too small. This can be seen immediately from  the fact that the antiderivative of the reciprocal is the global logarithm which is not definable in $\IR_\an$.
 \vs{0.2cm}
\hrule

\vs{0.4cm}
{\footnotesize{\itshape 2010 Mathematics Subject Classification:} 03C64, 03H05, 06F20, 12J25, 26E30, 28B15, 28E05, 32B20}
\newline
{\footnotesize{\itshape Keywords and phrases:} Real field with restricted analytic function, non-archimedean models, archimedean rank, globally subanalytic sets and functions, Lebesgue measure and integration}


However by the seminal work of Comte, Lion and Rolin [5, 17] (see also [11]), later extended by Cluckers and D. Miller [1, 2, 3], it is enough to add the logarithm to do integration on $\IR_\an$. 
Hence the existence of a logarithm is essential in this regard. A non-archimedean model $R$ of the theory $T_\an$ has a partial logarithm. 
But it cannot be in general extended to a reasonable global logarithm on $R$, for example if $R$ is a field of generalized power series (see Kuhlmann et al. [13]). 
There is the important construction of the field of transseries or logarithmic exponential series (see Van den Dries et al. [7]) where a global logarithm and inverse exponential on certain unions of an increasing sequence of fields of generalized power series fields are constructed (see also Kuhlmann [16] for another construction).\\
Since for integration only the logarithm has to be established our motivation is to extend the partial logarithm in a simple way as possible.
Assuming that the value group $\Gamma_R$ of $R$ with respect to the standard valuation induced by the ordering has finite archimedean rank $\ell$ (i.e. $\Gamma_R$ has $\ell$ many archimedean classes) we are able to define a global logarithm on $R$ with values in the polynomial ring $R[X]=R[X_1,\ldots,X_\ell]$ with excellent properties.
The construction is controlled by a so-called logarithmic datum, which is formulated in terms of $R$ and consists of a section for $R$ (i.e. a group embedding $\Gamma\hookrightarrow R_{>0}$) and of an embedding of $\Gamma_R$ into its Hahn group $\IR^\ell$ (see Prie\ss-Crampe [18, I \S 5]) where the latter is equipped with the antilexicographical ordering.
If the archimedean rank of the group $\Gamma_R$ equals its rational rank (i.e. its dimension as $\IQ$-vector space) then the logarithm is unique up to a unique isomorphism, hence indepent of the choice of the logarithmic datum.\\
A central result of the paper is that such a logarithm constructed as above can be embedded into a model of the theory $T_{\an,\exp}$ of the real field with restricted analytic functions and global exponential function (see again [6] for this theory). For this we use generalized fields of logarithmic-exponential series from [7]. With such embeddings we can define, after the choice of a logarithm, constructible (i.e. log-analytic) functions on $R$ with values in $R[X]$ in the sense of [1] with good lifting properties of constructible functions on the reals.\\
Having this in hand one can lift the results of [1, 2, 3, 5, 17] to the non-standard model $R$ to obtain a complete Lebesgue measure and integration theory for the globally subanalytic category with values in the polynomial ring $R[X]$. This extends and strengthen the construction in [14] where a Lebesgue theory was developed in the case of an archimedean value group. The non-archimedean Lebesgue measure and integration theory here has to be put in context with other integration theories in non-archimedean settings. Here the goal is to transfer the ubiquitious integration from analysis to non-archimedean orderings which is in its numerical character of quite different flavour than motivic integration in non-archimedean geometry (see for example Cluckers et al. [4] for the latter). But there are interesting similarities to recent work by Ducros et al. [8] describing non-archimedean integrals as limits of complex integrals since the Lebesgue integral here is also constructed via families of functions on the reals. They work in the highly saturated setting of Robinson's non-standard analysis whereas we want to extend the ground field as less as possible.\\
The paper is organized as follows. 
In the first preliminary section notations are introduced, necessary background on ordered abelian groups and the theories $T_\an$ and $T_{\an,\exp}$ is given and the general setting is established. In Section 2 the logarithmic functions are defined and analyzed. In Section 3 we prove how the construction can be embedded into models of $\IR_{\an,\exp}$.  
In Section 4 the concept of constructible functions is developed. In the final Section 5 we construct as an application a Lebesgue theory for the globally subanalytic category. We have tried on one hand to avoid overlaps with the presentation in [14] and on the other hand to elaborate the new features and to keep this paper self-contained. Section 5 ends with the discussion of the  strenghtening of the present construction compared to the one in [14] in the case of an archimedean value group.
  
\section{Preliminaries}

\subsection{Notations}

By $\IN=\big\{1,2,3,\ldots\big\}$ we denote the set of natural numbers and by $\IN_0=\big\{0,1,2,\ldots\big\}$ the set of natural numbers with $0$.\\
We set $\IR^*:=\{x\in \IR\mid x\neq 0\}, \IR_{>0}:=\{x\in \IR\mid x>0\}$ and $\IR_{\geq 0}:=\{x\in \IR\mid x\geq 0\}$. An element of $\IR_{>0}$ is called positive. For $a,b\in \IR$ with $a<b$ let
$]a,b[\;:=\{x\in \IR\mid a<x<b\}$ be the open and $[a,b]:=\{x\in \IR\mid a\leq x\leq b\}$ be the closed interval with endpoints $a$ and $b$.\\
Given a subset $A$ of $\IR^n$ we denote by $\mathbbm{1}_A$ the characteristic function of $A$.
For a function $f:\IR^n\to \IR$ we set $f_+:=\max(f,0)$ and $f_-:=\max(-f,0)$. 
For a function $f:\IR^m\times \IR^n\to \IR, (x,y)\mapsto f(x,y),$ and $y\in \IR^m$ we denote by $f(-,y)$ the function $\IR^m\to \IR, x\mapsto f(x,y)$.\\
Given $n\in \IN$ we denote by $M(n,\IR)$ the set of square matrizes with real entries and $n$ rows.\\
The same notions apply if $\IR$ is replaced by an arbitrary ordered ring.

\vs{0.2cm}
Given sets $U,V,X,Y$ a map $f:Y\to Y$ and a map $G:V^U\to Y^X$ we write $f\circ G$ for the map $V^U\to Y^X, h\mapsto f\circ G(h)$.

\vs{0.2cm}
Finally, by $\infty$ we denote an element that is bigger than every element of a given ordered set.

\newpage

\subsection{Ordered abelian goups and ordered rings}

Let $\Gamma$ be an ordered abelian group that is divisible (or equivalently, an ordered $\IQ$-vector space). Its dimension as a $\IQ$-vector space is called the {\bf rational rank} and is denoted by $\mathrm{rk}_\IQ(\Gamma)$.\\
Given $\gamma\in \Gamma$ we denote by $|\gamma|:=\max\{\gamma,-\gamma\}$ the {\bf absolute value} of $\gamma$.
We say that $\gamma,\delta\in \Gamma$ are {\bf archimedean equivalent} and write $a\sim b$ if there are $m,n\in\IN$ such that $|\gamma|\leq m|\delta|$ and $|\delta|\leq n|\gamma|$. 
Given $\gamma,\delta\in \Gamma_{>0}$ we set $\gamma\ll \delta$ if $n\gamma<\delta$ for all $n\in\IN$.\\
For $\gamma\in\Gamma\setminus \{0\}$ we denote by 
$[\gamma]:=\{\delta\in \Gamma\mid \gamma\sim \delta\}$ the archimedean class of $\gamma$. By $\Delta=\Delta(\Gamma):=\{[\gamma]\mid\gamma\in\Gamma\setminus\{0\}\}$ we denote the set of archimedean classes of $\Gamma$.
We equip $\Delta$ with the ordering induced by $\Gamma$. 
The cardinality of the set $\Delta$ is called the {\bf archimedean rank} of $\Gamma$ and is denoted by $\mathrm{rk}_\mathrm{arch}(\Gamma)$. 
Note that $\mathrm{rk}_\IQ(\Gamma)\geq\mathrm{rk}_\mathrm{arch}(\Gamma)$.
We have $\mathrm{rk}_\mathrm{arch}(\Gamma)=0$ if and only if $\Gamma=\{0\}$ and $\mathrm{rk}_\mathrm{arch}(\Gamma)=1$ if and only if $\Gamma\neq\{0\}$ is archimedean.
If $\Gamma$ has finite archimedean rank $\ell$ let $\delta_1<\ldots<\delta_\ell$ be the elements of $\Delta$.
We set 
$$w=w_\Gamma:\Gamma\to \big\{0,\ldots,\ell\big\},
x\mapsto \left\{\begin{array}{ccc}
k,&&x\neq 0\mbox{ and }[x]=\delta_k,\\
&\mbox{if}&\\
0,&&x=0.\end{array}\right.
$$

\vs{0.5cm}
{\bf 1.1 Example}

\vs{0.1cm}
Let $\ell\in\IN$. 
\begin{itemize}
\item[(1)] 
The group $\IQ^\ell$ equipped with the antilexicographical ordering is denoted by $\IQ_\mathrm{anlex}^\ell$. We have $\mathrm{rk}_\IQ(\IQ_\mathrm{anlex}^\ell)=\mathrm{rk}_\mathrm{arch}(\IQ_\mathrm{anlex}^\ell)=\ell$.
\item[(2)] 
The group $\IR^\ell$ equipped with the antilexicographical ordering is denoted by $\IR_\mathrm{anlex}^\ell$. We have $\mathrm{rk}_\IQ(\IR_\mathrm{anlex}^\ell)=\infty$ and $\mathrm{rk}_\mathrm{arch}(\IR_\mathrm{anlex}^\ell)=\ell$.
\end{itemize}
The archimedean classes of $\IQ_\mathrm{anlex}^\ell$ (respectively $\IR_\mathrm{anlex}^\ell$) are $[e_1]<\ldots<[e_\ell]$ where $e_k$ denotes the $k^\mathrm{th}$-unit vector for $k\in\{1,\ldots,\ell\}$.
For $x\in \IQ_\mathrm{anlex}^\ell$ (respectively $\IR_\mathrm{anlex}^\ell$) we have 
$$w(x)=\min\big\{k\in\{0,\ldots,\ell\}\;\big\vert\; x\in \IR e_1+\ldots+\IR e_k\big\}.$$
Given $x\neq 0$ we have $[x]=[e_{w(x)}]$.

\vs{0.2cm}
Moreover we set $\IR^0:=\{0\}$. 

\vs{0.5cm}
The following is a special case of the Hahn embedding theorem (see Prie\ss-Crampe [18, I \S 5]).

\vs{0.5cm}
\newpage
{\bf 1.2 Fact}

\vs{0.1cm}
{\it Assume that $\Gamma$ has finite archimedean rank $\ell$.
	\begin{itemize}
		\item[(1)] There is an order preserving group embedding $\Gamma\hookrightarrow \IR_\mathrm{anlex}^\ell$.
		\item[(2)] If $\tau,\tau':\Gamma\hookrightarrow  \IR_\mathrm{anlex}^\ell$ are order preserving group embeddings then there is an order isomorphism $\lambda:\IR_\mathrm{anlex}^\ell\stackrel{\cong}{\longrightarrow}\IR_\mathrm{anlex}^\ell$ such that
		$\tau'=\lambda\circ\tau$.
	\end{itemize}}

\vs{0.2cm}
Note that in the above situation an order preserving embedding $\tau:\Gamma\hookrightarrow \IR_\mathrm{anlex}^\ell$ is also valuation preserving, that is $w_\Gamma=w_{\IR_\mathrm{anlex}^\ell}\circ\tau$. 

\vs{0.5cm}
Let $A$ be a ring (commutative with unit) and $x=(x_1,\ldots,x_\ell),y=(y_1,\ldots,y_\ell)\in A^\ell$. We set $\big\langle x,y\big\rangle:=x_1y_1+\ldots+x_\ell y_\ell$.\\
Assume that $O$ is an ordered ring. Given $a,b\in O_{>0}$ we write $a\prec b$ if $a^m<b^n$ for all $m,n\in\IN$.\\
Assume that $R$ is an ordered field.
Let
$$\ma{O}_R:=\big\{a\in R\;\big\vert\;-n\leq a\leq n\mbox{ for some }n\in\IN\big\}$$
be the set of {\bf bounded} elements,
$$\mathfrak{m}_R:=\big\{a\in R\;\big\vert\;-1/n\leq a\leq 1/n\mbox{ for all }n\in\IN\big\}$$
the set of {\bf infinitesimal} elements
and
$$\ma{I}_R:=\big\{a\in R\;\big\vert\;a\geq n\mbox{ for all }n\in\IN\big\}$$
the set of {\bf infinitely large} elements of $R$, respectively.
Then $\ma{O}_R$ is a convex valuation ring of $R$ with maximal ideal
$\mathfrak{m}_R$.
Let $v_R:R^*\to R^*/\ma{O}_R^*$ denote the corresponding {\bf standard real valuation} with {\bf value group} $\Gamma_R:=R^*/\ma{O}_R^*$. Note that $\Gamma_R$ is divisible if $R$ is real closed.

\subsection{The real field with restricted analytic functions}

Let $\IR_\an$ be the real field with restricted analytic functions.
Its theory with respect to its natural language $\ma{L}_\an$ is denoted by $T_\an$.
Note that a model of $T_\an$ is in particular a real closed field with $\IR_\an$ as elementary substructure.
In Van den Dries et al. [6] the theory $T_\an$ has been deeply analyzed. We need the following facts.

\vs{0.5cm}
{\bf 1.3 Fact}

\vs{0.1cm}
{\it Let $\Gamma$ be an ordered abelian group that is divisible.
Then the power series field $\IR((t^\Gamma))$ can be made in a natural way into a model of $T_\an$.}

\vs{0.5cm}
Let $R$ be a model of $T_\an$. A section for $R$ is a group homomorphism $s:(\Gamma_R,+)\to (R_{>0},\cdot)$ such that $v_R\big(s(\gamma)\big)=\gamma$ for all $\gamma\in\Gamma_R$. Since $\Gamma_R$ is divisible there is a section for $R$. 

\vs{0.5cm}
{\bf 1.4 Fact}

\vs{0.1cm}
{\it Let $R$ be a model of $T_\an$ and let $s$ be a section for $R$. Then there is an $\ma{L}_\an$-embedding $\sigma:R\hookrightarrow \IR((t^{\Gamma_R}))$ such that $\sigma\big(s(\gamma)\big)=t^\gamma$ for all $\gamma\in\Gamma$.}

\vs{0.5cm}
Let $R$ be a model of $T_\an$. Let $n\in\IN$ and let $A$ be a subset of $R^n$. Then $A$ (respectively a function $f:A\to R$) is {\bf globally subanalytic} if it is $\ma{L}_\an$-definable.
By $\mathrm{Sub}_R(A)$ we denote the ring of the globally subanalytic functions on a globally subanalytic set $A$.\\
Let $R,\widetilde{R}$ be models of $T_\an$ such that $R$ is a substructure of $\widetilde{R}$. Then $R$ is an elementary substructure of $\widetilde{R}$. Given a globally subanalytic subset $A$ of $R^n$ or a globally subanalytic 
function $f:A\to R$ we denote by $A_{\widetilde{R}}$ and $f_{\widetilde{R}}:A_{\widetilde{R}}\to \widetilde{R}$ their canonical lifting to a globally subanalytic subset of $\widetilde{R}^n$ respectively a globally subanalytic function on $A_{\widetilde{R}}$.

\subsection{The real field with restricted analytic functions and exponentiation}

Let $\IR_{\an,\exp}$ be the real field with restricted analytic functions and exponentiation.
Its theory with respect to its natural language $\ma{L}_{\an,\exp}$ is denoted by $T_{\an,\exp}$.
In [6] the theory $T_{\an,\exp}$ has been deeply analyzed. 
Note that a model of $T_{\an,\exp}$ has $\IR_{\an,\exp}$ as an elementary substructure and is in particular a model of $T_\an$.

\vs{0.5cm}
{\bf 1.5 Proposition}

\vs{0.1cm}
{\it Let $R$ be a model of $T_\an$. Then there is a model $S$ of $T_{\an,\exp}$ extending $R$.}

\vs{0.1cm}
{\bf Proof:}
The fields $R$ and $\IR$ are both models of the complete theory $T_\an$.
By the Keisler-Shelah-Theorem (see Hodges [9, p. 245]) there is a set $I$ and an ultrafilter
$\ma{U}$ on $I$ such that $R^I/\ma{U}\cong \IR^I/\ma{U}$ as $\ma{L}_\an$-structures.
Clearly $S:=\IR^I/\ma{U}$ as an ultrapower of $\IR$ is a model of the $\IR_{\an,\exp}$-theory $T_{\an,\exp}$. Since $R$ is an elementary substructure of the ultrapower $R^I/\ma{U}$ we can view $R$ via the $\ma{L}_\an$-isomorphism $R^I/\ma{U}\cong S$ as a substructure of $S$.
\hfill$\Box$

\vs{0.5cm}
Let $S$ be a model of $T_{\an,\exp}$. Let $n\in\IN$ and let $A$ be a globally subanalytic subset of $S^n$. 
A function $f:A\to S$ is called {\bf constructible} if it is a finite sum of finite products of globally subanalytic functions and the logarithm of positive globally subanalytic functions on $A$ (compare with Cluckers and Dan Miller [1]).
By $\mathrm{Con}_S(A)$ we denote the ring of all constructible functions on a globally subanalytic set $A$.\\
Let $S,\widetilde{S}$ be models of $T_{\an,\exp}$ such that $S$ is a substructure of $\widetilde{S}$. Then $S$ is an elementary substructure of $\widetilde{S}$. Let $A$ be a globally subanalytic subset of $S^n$ and let $f:A\to S$ be a constructible function.
The canonical lifting $f_{\widetilde{S}}:A_{\widetilde{S}}\to \IR$ is a constructible function on $A_{\widetilde{S}}$.

\subsection{General assumption}

Throughout the rest of the paper we work under the following general assumption:

\begin{center}
	{\it Let $R$ be a model of $T_\an$ such that the value group $\Gamma=\Gamma_R$ has finite archimedean rank $\ell$.}
\end{center}

Note that $\ell=0$ if and only if $R=\IR$.\\
An example is given by the field of Puiseux series over $\IR$ in $\ell$ variables.

\vs{0.5cm}
{\bf 1.6 Example}

By 
$$\mathbb{P}_\ell:=\Big\{t_1^{-k_1/p}\cdot\ldots\cdot t_\ell^{-k_\ell/p}f(t_1^{1/p},\ldots,t_\ell^{1/p})\;\Big\vert\;f\in \IR[[t_1,\ldots,t_\ell]], k_1,\ldots,k_\ell\in \IN_0\mbox{ and }p\in \IN\Big\}$$
we denote the field of (formal) {\bf Puiseux serie}s over $\IR$ in $\ell$ variables. 
We order $\mathbb{P}_\ell$ by letting 
$$0< t_\ell\prec \ldots\prec t_1\prec 1.$$
Then $\mathbb{P}_\ell$ is a non-archimedean model of $T_\an$ with value group $\Gamma_{\mathbb{P}_\ell}=\IQ_\mathrm{anlex}^\ell$. We have that $v_{\mathbb{P}_\ell}(t_k)=e_k$ for $k\in \{1,\ldots,\ell\}$.

\vs{0.5cm}
We denote by $R[X]=R[X_1,\ldots,X_\ell]$ the polynomial ring over $R$ in $\ell$ variables.
We order $R[X]$ by setting
$$1\prec X_1\prec\ldots\prec X_\ell \prec\ma{I}_R.$$

\section{Logarithms}

The $T_\an$-model $R$ carries a {\bf partial logarithm}
$$\log:\IR_{>0}+\mathfrak{m}_R\to R, a+m\mapsto \log(a)+L(m/a),$$
where $L$ denotes the logarithmic series.

\vs{0.5cm}
{\bf 2.1 Remark}

\vs{0.1cm}
{\it The partial logarithm is an isomorphism $\big(\IR_{>0}+\mathfrak{m}_R,\cdot\big)\stackrel{\cong}{\longrightarrow}\big(\ma{O}_R,+\big)$ of ordered groups.
Its inverse is the partial exponential map 
$$\exp:\ma{O}_R\to \IR_{>0}+\mathfrak{m}_R, a+m\mapsto e^aE(m),$$ 
where $a\in \IR,m\in \mathfrak{m}_R$, and $E$ denotes the exponential series.}

\vs{0.5cm}
{\bf 2.2 Definition}

\vs{0.1cm}
A {\bf logarithmic datum} for $R$ is a tuple $\big(s,\tau\big)$ where $s:\Gamma\hookrightarrow R_{>0}$ is a section for $R$ and $\tau:\Gamma\hookrightarrow \IR_\mathrm{anlex}^\ell$ is an order preserving embedding.

\vs{0.5cm}
{\bf 2.3 Remark}

\vs{0.1cm}
There is a logarithmic datum for $R$.

\vs{0.1cm}
{\bf Proof:}

\vs{0.1cm}
In Section 1.3 we have noted that there is a section $s$ for $R$. Note that there is an order preserving embedding $\tau:\Gamma\hookrightarrow \IR_\mathrm{anlex}^\ell$ has been formulated in Fact 1.2(1).

\hfill$\Box$

\vs{0.5cm}
{\bf 2.4 Definition}

\vs{0.1cm}
Let $\mu=\big(s,\tau\big)$ be a logarithmic datum for $R$. The $\mu$-logarithm $\log_\mu: R_{>0}\to R[X]$ is defined as follows.
Let $x\in R_{>0}$. Then $x/\big(s(v_R(x))\big)\in \IR_{>0}+\mathfrak{m}_R$. 
We set  
$$\log_\mu(x)=-\Big\langle \tau(v_R(x)),X\Big\rangle +\log\Big(\frac{x}{s(v_R(x))}\Big).$$

\vs{0.5cm}
{\bf 2.5 Proposition}

\vs{0.1cm}
{\it Let $\mu=\big(s,\tau\big)$ be a logarithmic datum for $R$. The following holds.
\begin{itemize}
	\item[(1)] The $\mu$-logarithm extends the partial logarithm on $R$.
	\item[(2)] The image of the $\mu$-logarithm is $\big\langle\tau(\Gamma),X\big\rangle+\ma{O}_R$. 
\end{itemize}}

\vs{0.1cm}
{\bf Proof:}

\vs{0.1cm}
(1): Let $x\in \IR_{>0}+\mathfrak{m}_R$. Then $v_R(x)=0$ and hence $\tau(v_R(x))=0$. Since $s(0)=1$ we obtain by the above definition that $\log_\mu(x)=\log(x)$.

\vs{0.2cm}
(2): That $\log_\mu(R_{>0})\subset \big\langle\tau(\Gamma),X\big\rangle+\ma{O}_R$ is clear by the definition of the $\mu$-logarithm. 
Let $\gamma\in\Gamma$ and let $g\in \ma{O}_R$. By Remark 2.1 there is $u\in \IR_{>0}+\mathfrak{m}_R$ such that $\log(u)=g$. Let $x:=s(-\gamma)u$. Then $\log_\mu(x)=\big\langle \tau(\gamma),X\big\rangle+g$.
\hfill$\Box$

\vs{0.5cm}
Note that the logarithm induced by a logarithmic datum is surjective only in the case $\ell=0$; i.e. in the case of the real field. 

\vs{0.5cm}
{\bf 2.6 Theorem}

\vs{0.1cm}
{\it Let $\mu=\big(s,\tau\big)$ be a logarithmic datum for $R$. The $\mu$-logarithm is an order preserving group embedding $\big(R_{>0},\cdot\big)\hookrightarrow \big(R[X],+\big)$.}

\vs{0.1cm}
{\bf Proof:}

\vs{0.1cm}
To show that $\log_\mu$ is a homomorphism let $x,x'\in R_{>0}$.
Let $\gamma:=v_R(x)$ and $\gamma':=v_R(x')$. Then $v_R(xx')=\gamma+\gamma'$ and $s(\gamma+\gamma')=s(\gamma)s(\gamma')$.
We obtain
\begin{eqnarray*}
	\log_\mu(xx')&=&-\Big\langle \tau(\gamma+\gamma'),X\Big\rangle+\log\Big(\frac{xx'}{s(\gamma)s(\gamma')}\Big)\\
	&\stackrel{\mathrm{Remark\,2.1}}{=}&-\Big\langle \tau(\gamma),X\Big\rangle +\log\Big(\frac{x}{s(\gamma)}\Big)-\Big\langle \tau(\gamma'),X\Big\rangle+\log\Big(\frac{x'}{s(\gamma')}\Big)\\
	&=&\log_\mu(x)+\log_\mu(x').
\end{eqnarray*}
To conclude that $\log_\mu$ is an order preserving embedding it is enough to show the following. Let $x\in R$ with $x>1$. Then $\log_\mu(x)>0$.
Set $\gamma:=v_R(x)$. Then $\gamma\leq 0$. If $\gamma=0$ we have $x\in \IR_{>0}+\mathfrak{m}_R$. Hence we are done by Proposition 2.5(1) and Remark 2.1.
In the case that $\gamma<0$ we have by the definition of $\log_\mu$ that
$\log_\mu(x)=\big\langle c ,X\big\rangle +w$ where $c:=-\tau(\gamma)\in \IR^\ell_\mathrm{anlex}$ is positive and $w\in \ma{O}_R$.
By the definition of the ordering on $R[X]$ we have that $\big\langle c,X\big\rangle >\ma{O}_R$. Hence $\log_\mu(x)>1$.
\hfill$\Box$

\vs{0.5cm}
{\bf 2.7 Remark}

\vs{0.1cm}
In the case $R=\IR$ (i.e. $\ell=0$) there is exactly one logarithmic datum. The corresponding logarithm is the usual logarithm.

\vs{0.5cm}
{\bf 2.8 Example}

\vs{0.1cm}
In the case of the field $\mathbb{P}_\ell$ of Puiseux series in $\ell$ variables let $e$ be the section $\IQ_\mathrm{anlex}^\ell\hookrightarrow \mathbb{P}_{\ell,>0}, \gamma\mapsto t^\gamma,$ and let $\iota:\IQ_\mathrm{anlex}^\ell\hookrightarrow \IR_\mathrm{anlex}^\ell$ be the inclusion.
The tuple $\nu=\nu_\ell=(e,\iota)$ is called the canonical logarithmic datum for $\mathbb{P}_\ell$.
We have that $\log_\nu(t_k^{-1})=X_k$ for $k\in\{1,\ldots,\ell\}$.

\vs{0.5cm}
We compare distinct logarithmic data.

\vs{0.5cm}
{\bf 2.9 Definition}

\vs{0.1cm}
Two logarithmic data $\mu,\mu'$ for $R$ are called {\bf equivalent} if there is a ring homomorphism $\Phi:R[X]\to R[X]$ such that $\Phi\circ\log_\mu=\log_{\mu'}$.

\vs{0.5cm}
{\bf 2.10 Proposition}

\vs{0.1cm}
{\it Let $\mu,\mu'$ be logarithmic data for $R$ which are equivalent. Let $\Phi:R[X]\stackrel{\cong}{\longrightarrow}R[X]$ be a ring homomorphism such that $\Phi\circ\log_\mu=\log_{\mu'}$.
	The following holds.
	\begin{itemize}
		\item[(1)] The restriction of $\Phi$ to $R$ is the identity.
		\item[(2)] There are a strictly lower triangular matrix $A\in M(\ell,\IR)$, a diagonal matrix $D\in \mathrm{M}(\ell,\IR)$ with positive elements in the diagonal and $g\in \ma{O}_R^\ell$ such that
		$\Phi(X)=(A+D)X+g$.
		\item[(3)] The homomorphism $\Phi$ is an order isomorphism. 
\end{itemize}}
{\bf Proof:}

\vs{0.1cm}
(1): By Proposition 2.5(1) and Remark 2.1 we have that $\Phi|_{\ma{O}_R}=\mathrm{id}_{\ma{O}_R}$. Since $\Phi$ is a ring homomorphism we get that $\Phi|_R=\mathrm{id}_R$. 

\vs{0.2cm}
(2): Let $\mu=\big(s,\tau\big)$ and $\mu'=\big(s',\tau'\big)$. Let $k\in\{1,\ldots,\ell\}$. We choose $\gamma_k\in \Gamma_{<0}$ such that $w(\tau(\gamma_k))=k$. Then $w(\tau'(\gamma_k))=k$. Let $x_k:=s(\gamma_k)\in R_{>0}$.
We have $\log_\mu(x_k)=-\big\langle \tau(\gamma_k),X\big\rangle$ and $ \log_{\mu'}(x_k)=-\big\langle\tau'(\gamma_k),X\big\rangle+h_k$ for some $h_k\in\ma{O}_R$. 
Hence
$$-\Big\langle \tau(\gamma_k),\Phi(X)\Big\rangle=
\Phi\Big(-\Big\langle\tau(\gamma_k),X\Big\rangle\Big)=
-\Big\langle\tau'(\gamma_k),X\Big\rangle +h_k.$$
Since $\gamma_k<0$ and $w(\tau(\gamma_k))=w(\tau'(\gamma_k))=k$ we find
$a_{k1},\ldots,a_{kk}\in \IR$ with $a_{kk}>0$ and $b_{k1},\ldots,b_{kk}\in \IR$ with $b_{kk}>0$ such that 
$-\tau(\gamma_k)=a_{k1}e_1+\ldots+a_{kk}e_k$ and $-\tau'(\gamma_k)=b_{k1}e_1+\ldots+b_{kk}e_k$.
Hence 
$$a_{k1}\Phi(X_1)+\ldots+a_{kk}\Phi(X_k)=b_{k1}X_1+\ldots+b_{kk}X_k+h_k.$$
Inductively we get that for every $k\in \{1,\ldots,\ell\}$ there are
$c_{k1},\ldots,c_{kk}\in \IR$ with $c_{kk}>0$ and $g_k\in \ma{O}_R$ such that
$$\Phi(X_k)=c_{k1}X_1+\ldots+c_{kk}X_k+g_k.$$

\vs{0.2cm}
(3): The fact that $\Phi$ is an isomorphism follows from (1) and (2). By the ordering given on $R[X]$ we obtain that $\Phi$ is order preserving.
\hfill$\Box$

\vs{0.5cm}
{\bf 2.11 Corollary}

\vs{0.1cm}
{\it Let $\mu,\mu'$ be logarithmic data for $R$ which are equivalent. Then there is exactly one ring homomorphism (resp. isomorphism) $\Phi:R[X]\to R[X]$ such that $\Phi\circ\log_\mu=\log_{\mu'}$. }

\vs{0.1cm}
{\bf Proof:}

\vs{0.1cm}
This follows from Proposition 2.10(1) and the proof of Proposition 2.10(2).
\hfill$\Box$

\vs{0.5cm}
{\bf 2.12 Definition}

\vs{0.1cm}
Let $\mu,\mu'$ be logarithmic data for $R$ which are equivalent. The unique ring homomorphism $\Phi:R[X]\stackrel{\cong}{\longrightarrow}R[X]$ such that $\Phi\circ\log_\mu=\log_{\mu'}$ is called
the {\bf logarithmic connection} between $\mu$ and $\mu'$.

\vs{0.5cm}
{\bf 2.13 Corollary}

\vs{0.1cm}
{\it The equivalence of logarithmic data is an equivalence relation.} 

\vs{0.1cm}
{\bf Proof:}

\vs{0.1cm}
This follows from Proposition 2.10(3).
\hfill$\Box$

\vs{0.5cm}
\newpage
{\bf 2.14 Proposition}

\vs{0.1cm}
{\it Let $s:\Gamma\hookrightarrow R_{>0}$ be a section for $R$ and let $\tau,\tau':\Gamma\hookrightarrow \IR_\mathrm{anlex}^\ell$ be order preserving.
	Then the logarithmic data $\mu:=\big(s,\tau\big)$ and $\mu'=\big(s,\tau'\big)$ are equivalent.}

\vs{0.1cm}
{\bf Proof:}

\vs{0.1cm}
By Fact 1.2(2) there is an order isomorphism $\lambda:\IR_\mathrm{anlex}^\ell\stackrel{\cong}{\longrightarrow}\IR_\mathrm{anlex}^\ell$ such that $\tau'=\lambda\circ\tau$.
Let $\Phi:R[X]\stackrel{\cong}{\longrightarrow}R[X]$ be the isomorphism extending the identity on $R$ such that $X_k\mapsto \big\langle\lambda(e_k),X\big\rangle$ for $k\in\{1,\ldots,\ell\}$. 
We show that $\Phi$ is the logarithmic connection between $\mu$ and $\mu'$.
Let $x\in R_{>0}$ and set $\gamma:=v_R(x)$. Then
$$\Phi\big(\log_\mu(x)\big)=\Phi\Big(-\Big\langle \tau(\gamma),X\Big\rangle+\log\Big(\frac{x}{s(\gamma)}\Big)\Big)
=-\Big\langle \tau(\gamma),\Phi(X)\Big\rangle+\log\Big(\frac{x}{s(\gamma)}\Big).$$
Since 
$$\log_{\mu'}(x)=-\big\langle\tau'(\gamma),X\big\rangle+\log\Big(\frac{x}{s(\gamma)}\Big)$$ 
we have to show that
$\big\langle\tau(\gamma),\Phi(X)\big\rangle=\big\langle\tau'(\gamma),X\big\rangle$.
We compute
\begin{eqnarray*}
	\Big\langle\tau(\gamma),\Phi(X)\Big\rangle	
	&=&\sum_{k=1}^\ell\tau_k(\gamma)\Phi(X_k)
	=\sum_{k=1}^\ell\tau_k(\gamma)\Big\langle \lambda(e_k),X\Big\rangle\\
	&=&\sum_{k=1}^\ell\tau_k(\gamma)\Big(\sum_{j=1}^\ell\lambda_j(e_k)X_j\Big)
	=\sum_{j=1}^\ell\Big(\sum_{k=1}^\ell\lambda_j(e_k)\tau_k(\gamma)\Big)X_j\\
	&=&\sum_{j=1}^\ell\lambda_j(\tau(\gamma))X_j
	=\Big\langle \lambda\big(\tau(\gamma)\big),X\Big\rangle\\
	&=&\Big\langle \tau'(\gamma),X\Big\rangle.\\
\end{eqnarray*}
\hfill$\Box$

\vs{0.5cm}
{\bf 2.15 Definition}

\vs{0.1cm}
We say that $R$ has a {\bf unique logarithm} if any two logarithmic data for $R$ are equivalent. 

\vs{0.5cm}
{\bf 2.16 Theorem}

\vs{0.1cm}
{\it The following are equivalent.
	\begin{itemize}
		\item[(i)] $\mathrm{rk}_\IQ(\Gamma)=\mathrm{rk}_\mathrm{arch}(\Gamma)$.
		\item[(ii)] $R$ has a unique logarithm.
\end{itemize}}

\vs{0.1cm}
{\bf Proof:}

\vs{0.1cm}
(i) $\Longrightarrow$ (ii):
Let $\mu=\big(s,\tau\big)$ and $\mu'=\big(s',\tau'\big)$ be logarithmic data for $R$. We have to show that $\mu$ and $\mu'$ are eqivalent. By Proposition 2.14 and Corollary 2.13 we can assume that $\tau=\tau'$.
Choose $\gamma_1,\ldots,\gamma_\ell\in\Gamma$ such that $0<\gamma_1\ll \ldots\ll \gamma_\ell$. 
Since $\mathrm{rk}_\IQ(\Gamma)=\mathrm{rk}_\mathrm{arch}(\Gamma)=\ell$ we have that $\{\gamma_1,\ldots,\gamma_\ell\}$ is a $\IQ$-basis of $\Gamma$.
For $k\in\{1,\ldots,\ell\}$ let 
$$h_k:=\log\Big(\frac{s'(\gamma_k)}{s(\gamma_k)}\Big)\in \ma{O}_R$$ 
and set $h:=(h_1,\ldots,h_\ell)\in\ma{O}_R^\ell$. 
We have that $w(\tau(\gamma_k))=k$ for every $k\in\{1,\ldots,\ell\}$. Hence the lower triangular matrix 
$$A:=\left(\begin{array}{c}
\tau(\gamma_1)\\
\cdot\\
\cdot\\
\cdot\\
\tau(\gamma_\ell)\\
\end{array}\right)\in M(\ell,\IR)$$ is invertible.
Let $g:=A^{-1}h\in \ma{O}_R^\ell$.  
Let $\Phi:R[X]\stackrel{\cong}{\longrightarrow}R[X]$ be the isomorphism extending the identity on $R$ such that $X_k\mapsto X_k+g_k$ for every $k\in\{1,\ldots,\ell\}$. 
We show that $\Phi$ is the logarithmic connection between $\mu$ and $\mu'$. 
Let $x\in R_{>0}$ and set $\gamma:=v_R(x)$. 
We have
\begin{eqnarray*}
	\Phi\Big(\log_\mu(x)\Big)&=&\Phi\Big(-\Big\langle\tau(\gamma),X\Big\rangle+\log\Big(\frac{x}{s(\gamma)}\Big)\Big)\\
	&=&-\Big\langle\tau(\gamma),\Phi(X)\Big\rangle +\log\Big(\frac{x}{s(\gamma)}\Big)\\
	&=&-\Big\langle\tau(\gamma),X\Big\rangle-\Big\langle\tau(\gamma),g\Big\rangle+\log\Big(\frac{x}{s(\gamma)}\Big)
\end{eqnarray*}
and 
$$\log_{\mu'}(x)=-\Big\langle\tau(\gamma),X\Big\rangle +\log\Big(\frac{x}{s'(\gamma)}\Big).$$
Hence it is enough to show that 
$$\log\Big(\frac{x}{s'(\gamma)}\Big)=\log\Big(\frac{x}{s(\gamma)}\Big)-\Big\langle\tau(\gamma),g\Big\rangle.$$ 
Applying the logarithmic rule this is equivalent to  
$$\log\Big(\frac{s'(\gamma)}{s(\gamma)}\Big)=\Big\langle\tau(\gamma),g\Big\rangle.$$
There are $q_1,\ldots,q_\ell\in \IQ$ such that $\gamma=q_1\gamma_1+\ldots+q_\ell\gamma_\ell$.
We have that
\begin{eqnarray*}
	\log\Big(\frac{s'(\gamma)}{s(\gamma)}\Big)&=&\log\Big(\frac{s'(q_1\gamma_1+\ldots+q_\ell\gamma_\ell)}{s(q_1\gamma_1+\ldots+q_\ell\gamma_\ell)}\Big)\\
	&=&\log\Big(\frac{s'(\gamma_1)^{q_1}\cdot\ldots\cdot s'(\gamma_\ell)^{q_\ell}}{s(\gamma_1)^{q_1}\cdot\ldots\cdot s(\gamma_\ell)^{q_\ell}}\Big)\\
	&=&q_1\log\Big(\frac{s'(\gamma_1)}{s(\gamma_1)}\Big)+\ldots+q_\ell\log\Big(\frac{s'(\gamma_\ell)}{s(\gamma_\ell)}\Big)\\
	&=&q_1h_1+\ldots+q_\ell h_\ell=\Big\langle q,h\Big\rangle\\
\end{eqnarray*}
and that
\begin{eqnarray*}
	\Big\langle\tau(\gamma),g\Big\rangle&=&\Big\langle\tau(q_1\gamma_1+\ldots+q_\ell\gamma_\ell),g\Big\rangle\\
	&=&q_1\Big\langle\tau(\gamma_1),g\Big\rangle+\ldots+q_\ell\Big\langle\tau(\gamma_\ell),g\Big\rangle\\
	&=&\Big\langle q,Ag\Big\rangle=\Big\langle q,h\Big\rangle.\\
\end{eqnarray*}
Hence the equality is established.

\vs{0.2cm}
(ii) $\Longrightarrow$ (i):
Assume that $\mathrm{rk}_\IQ(\Gamma)>\mathrm{rk}_\mathrm{arch}(\Gamma)=\ell$. 
Then there are $k\in\{1,\ldots,\ell\}$ and 
there are $\gamma,\delta$ in $\Gamma_{<0}$ such that $\gamma$ and $\delta$ are $\IQ$-linearly independent but $w(\gamma)=w(\delta)=k$.
Let $a,b\in R_{>0}$ such that $v_R(a)=\gamma$ and $v_R(b)=\delta$.
There is a section $s:\Gamma\hookrightarrow R_{>0}$ such that $s(\gamma)=a$ and $s(\delta)=b$ and there is a section $s':\Gamma\hookrightarrow R_{>0}$ such that $s'(\gamma)=a$ and $s'(\delta)=b/2$.
We find an order preserving embedding 
$\tau:\Gamma\hookrightarrow\IR_\mathrm{anlex}^\ell$ such that $\tau(\gamma)=-e_k$ and $\tau(\delta)=-\zeta e_k$ for some irrational $\zeta\in \IR_{>0}$. We show that the logarithmic data $\mu:=\big(s,\tau\big)$ and $\mu':=\big(s',\tau\big)$ are not equivalent.
We have 
$$\log_\mu(a)=X_k \mbox{ and }\log_{\mu'}(a)=X_k,$$
$$\log_\mu(b)=\zeta X_k\mbox{ and }\log_{\mu'}(b)=\log_{\mu'}\big(2s'(\delta)\big)=\zeta X_k+\log 2.$$
Hence if there were a homomorphism $\Phi:R[X]\to R[X]$ which fixes $R$ and fulfils $\Phi\circ\log_\mu=\log_{\mu'}$ we would obtain by the first line that $\Phi(X_k)=X_k$ and by the second that $\Phi(X_k)=X_k+\zeta^{-1}\log 2$, contradiction.
\hfill$\Box$

\newpage
\section{Embedding logarithms into models of $T_{\an,\exp}$}

We want to embed the logarithm induced by a logarithmic datum into a model of $T_{\an,\exp}$. Given a model $S$ of $T_{\an,\exp}$ we denote its logarithm by $\log_S$.

\vs{0.5cm}
{\bf 3.1 Definition}

\vs{0.1cm}
Let $\mu$ be a logarithmic datum for $R$. A {\bf $T_{\an,\exp}$-companion} for $\mu$ is a tuple $\big(S,\rho,Y\big)$ where $S$ is a model of $T_{\an,\exp}$,
$\rho:R\hookrightarrow S$ is an $\ma{L}_\an$-embedding and $Y=(Y_1,\ldots,Y_\ell)\in S^\ell$ such that the following conditions hold.
\begin{itemize}
	\item[(a)] The homomorphism $\rho^Y :R[X]\to S$ which extends $\rho$ and maps $X$ to $Y$ is an order preserving embedding.
	\item[(b)] The diagram 
\begin{center}
	
\begin{xy}\;\;\;\;\;\;\;\;\;\;\;\;\;\;\;\;\;\;\;\;\;\;\;\;\;\;\;\;\;\;\;\;\;\;\;\;\;\;\;\;\;\;\;\;\;\;\;\;\;\;\;\;\;\;\;\;\;\;\;\;
	\xymatrix{
	R_{>0}\ar[r]^{\log_\mu}\ar[d]_\rho& R[X]\ar[d]^{\rho^Y}\\
	S_{>0}\ar[r]^{\log_S}& S
}
\end{xy}
\end{center}
commutes.
\end{itemize}

\vs{0.2cm}
{\bf 3.2 Theorem}

\vs{0.1cm}
{\it A logarithmic datum has a $T_{\an,\exp}$-companion.} 

\vs{0.1cm}
{\bf Proof:}

\vs{0.1cm}
Let $\mu=\big(s,\tau\big)$ be a logarithmic datum for $R$. 

\vs{0.2cm}
{\bf Step 1:}
We set $Z:=(t_{\ell-1},\ldots,t_1,1)\in\mathbb{P}_{\ell-1}$ and define
$$\kappa=\kappa_\tau:\Gamma\to \mathbb{P}_{\ell-1}, \gamma\mapsto \Big\langle \tau(\gamma),Z\Big\rangle.$$
Then $\kappa$ is an order preserving embedding $\Gamma\hookrightarrow (\mathbb{P}_{\ell-1},+)$.

\vs{0.2cm}
{\bf Step 2:}
By Proposition 1.5 we find a model $K$ of $T_{\an,\exp}$ extending $\mathbb{P}_{\ell-1}$. 
Let $S$ be the field $K((t))^{\mathrm{LE}}$ of logarithmic-exponential series over $K$ from [7] (where $t$ denotes a new unknown different from $t_1,\ldots,t_{\ell-1})$.
Then $S$ is a model of $T_{\an,\exp}$.
It contains the $T_\an$-model $K((t^K))$ (cf. Fact 1.3).\\
We set  
$$\lambda:\IR((t^\Gamma))\hookrightarrow K((t^K))\subset S,
\sum_{\gamma\in \Gamma}a_{\gamma}t^\gamma\mapsto \sum_{\gamma\in\Gamma}a_\gamma t^{\kappa(\gamma)}\in K((t^K)).$$
Then $\lambda:\IR((t^\Gamma))\hookrightarrow S$ is an $\ma{L}_\an$-embedding.

\vs{0.2cm}
{\bf Step 3:}
By Fact 1.4 there is an $\ma{L}_\an$-embedding $\sigma:R\hookrightarrow\IR((t^\Gamma))$ such that $\sigma(s(\gamma))=t^\gamma$ for all $\gamma\in\Gamma$.
Let $\rho:=\lambda\circ \sigma$. We choose
$Y:=\log_S(t^{-1})Z\in S^\ell$ and show that $\big(S,\rho,Y\big)$ is a $T_{\an,\exp}$-companion for $\mu$.\\
Since $K_{>0}\prec\log_S(t^{-1})\prec t^{K_{<0}}$ we get that $\rho^Y:R[X]\to S$ is an order preserving embedding.
Let $x\in R_{>0}$. We have to show that $\log_S\big(\rho(x)\big)=\rho^Y\big(\log_\mu(x)\big)$. Let $\gamma:=v_R(x)$ and let $y:=x/s(\gamma)$. Then there are $a\in \IR_{>0}$ and $h\in \mathfrak{m}_R$ such that $y=a(1+h)$. 
Denoting the logarithmic series by $L$ we obtain by the construction of $K((t))^{\mathrm{LE}}$ and the fact that $\rho$  is an $\ma{L}_\an$-embedding that 
\begin{eqnarray*}
\log_S\big(\rho(x)\big)&=&\log_S\big(\lambda(\sigma(x))\big)=\log_S\big(\lambda(at^\gamma (1+\sigma(h)))\big)\\
&=&\log_S\big(at^{\kappa(\gamma)}(1+\rho(h))\big)\\
&=&\log(a)+\log_S(t^{\kappa(\gamma)})+L\big(\rho(h)\big)\\
&=&\log(a)-\kappa(\gamma)\log_S(t^{-1})+\rho\big(L(h)\big)\\
&=&\log(a)-\big\langle\tau(\gamma),Z\big\rangle \log_S(t^{-1})+\rho\big(L(h)\big)\\
&=&\log(a)-\big\langle\tau(\gamma),Y\big\rangle +\rho\big(L(h)\big)\\
&=&-\big\langle\tau(\gamma),Y\big\rangle +\rho\big(\log(y)\big)\\
&=&\rho^Y\Big(-\big\langle\tau(\gamma),X\big\rangle+\log(y)\Big)\\
&=&\rho^Y\big(\log_\mu(x)\big).
\end{eqnarray*}
\hfill$\Box$

\section{Constructible functions}

Let $\mu=\big(s,\tau\big)$ be a logarithmic datum for $R$.

\vs{0.5cm}
{\bf 4.1 Definition}

\vs{0.1cm}
Let $A$ be a globally subanalytic subset of $R^n$. A function $f:A\to R[X]$ is called {\bf $\mathbf{\mu}$-constructible}  
if it is a finite sum of finite products of globally subanalytic functions and the $\mu$-logarithm of positive globally subanalytic functions on $A$.
We denote by $\mathrm{Con}_{R,\mu}(A)$ the set of $\mu$-constructible functions on $A$.

\vs{0.5cm}
{\bf 4.2 Remark}

\vs{0.1cm}
Let $A$ be a globally subanalytic subset of $R^n$. Then $\mathrm{Con}_{R,\mu}(A)$ is an $R$-algebra extending $\mathrm{Sub}_R(A)$.

\vs{0.5cm}
Let $\Sigma:=\big(S,\rho,Y\big)$ be a $T_{\an,\exp}$-companion for $\mu$.
We have the following commuting diagram:

\begin{center}
	\begin{xy}\;\;\;\;\;\;\;\;\;\;\;\;\;\;\;\;\;\;\;\;\;\;\;\;\;\;\;\;\;\;\;\;\;\;\;\;\;\;\;\;\;\;\;\;\;\;\;\;\;\;\;\;\;\;\;\;\;\;\;\;
		\xymatrix{
			R\ar[rr]^{\rho}& &S\\
			&\IR.\ar[lu]^{\mathrm{incl}}\ar[ru]_{\mathrm{incl}}&
		}
	\end{xy}
\end{center}

We consider $S$ via $\rho$ as an $R$-algebra.
Since $\rho$ is an $\ma{L}_\an$-embedding we have that $\rho(R)$ is a $T_\an$-model.
Let $n\in\IN$. We denote the map $R^n\to S^n$ given componentwise by $\rho$ again by $\rho$. Let $A$ be a globally subanalytic subset of $R^n$.
We have a canonical isomorphism
$$\mathrm{Sub}_R(A)\hookrightarrow \mathrm{Sub}_{\rho(R)}(\rho(A)), f\mapsto \rho\circ f\circ\rho^{-1},$$
and a canonical lifting
$$\mathrm{Sub}_{\rho(R)}(\rho(A))\hookrightarrow \mathrm{Sub}_S(\rho(A)_S), g\mapsto g_S,$$
(cf. Section 1.3).
We write $A_\Sigma$ for $\rho(A)_S$ and denote the composition of these two maps by
$$\mathfrak{E}_\Sigma^A:\mathrm{Sub}_R(A)\hookrightarrow \mathrm{Sub}_S(A_\Sigma),f\mapsto f_\Sigma.$$
Note that $\mathfrak{E}_\Sigma^A$ is an embedding of $R$-algebras.
By construction, the following diagram commutes for $f\in \mathrm{Sub}_R(A)$:
\begin{center}
	\begin{xy}\;\;\;\;\;\;\;\;\;\;\;\;\;\;\;\;\;\;\;\;\;\;\;\;\;\;\;\;\;\;\;\;\;\;\;\;\;\;\;\;\;\;\;\;\;\;\;\;\;\;\;\;\;\;\;\;\;\;\;\;
		\xymatrix{
			A\ar[r]^{f}\ar[d]_{\rho}& R\ar[d]^{\rho}\\
			A_\Sigma\ar[r]^{f_\Sigma}& S.
		}
	\end{xy}
\end{center}

In view of the definition of an $T_{\an,\exp}$-companion we extend the embedding in a canonical way to $\mu$-constructible functions.

\vs{0.5cm}
{\bf 4.3 Construction}

The map $$\mathfrak{E}_{\mu,\Sigma}^A:\mathrm{Con}_{R,\mu}(A)\hookrightarrow \mathrm{Con}_S(A_\Sigma), f\mapsto f_\Sigma,$$ 
is defined as follows.
Let $f\in\mathrm{Con}_{R,\mu}(A)$.
There are $p,q\in \IN_0$ and for $i\in \{1,\ldots,p\},j\in\{1,\ldots,q\}$ there are $\varphi_i\in \mathrm{Sub}_R(A)$ and positive $\psi_{ij}\in\mathrm{Sub}_R(A)$ such that
$$f=\sum_{i=1}^p\varphi_i\Big(\prod_{j=1}^q\log_\mu\big(\psi_{ij}\big)\Big).$$
We set
$$f_\Sigma:=\sum_{i=1}^p\big(\varphi_i)_\Sigma\Big(\prod_{j=1}^q\log_S\big(\big(\psi_{ij}\big)_\Sigma\big)\big)\Big).$$

\vs{0.1cm}
{\bf Proof of well-definition:}

\vs{0.1cm}
This can be shown as in [14, Proposition 4.17 \& Theorem 4.18].
\hfill$\Box$

\vs{0.5cm}
\newpage
{\bf 4.4 Remark}

\vs{0.1cm}
The map $\mathfrak{E}_{\mu,\Sigma}^A:\mathrm{Con}_{R,\mu}(A)\hookrightarrow \mathrm{Con}_S(A_\Sigma)$ is an $R$-algebra homomorphism that extends $\mathfrak{E}_\Sigma^A$. 

\vs{0.5cm}
{\bf 4.5 Proposition}

\vs{0.1cm}
{\it The following diagram commutes for $f\in \mathrm{Con}_{R,\mu}(A)$:}
\begin{center}
	\begin{xy}\;\;\;\;\;\;\;\;\;\;\;\;\;\;\;\;\;\;\;\;\;\;\;\;\;\;\;\;\;\;\;\;\;\;\;\;\;\;\;\;\;\;\;\;\;\;\;\;\;\;\;\;\;\;\;\;\;\;\;\;
		\xymatrix{
			A\ar[r]^{f}\ar[d]_{\rho}& R[X]\ar[d]^{\rho^Y}\\
			A_\Sigma\ar[r]^{f_\Sigma}& S.
		}
	\end{xy}
\end{center}

\vs{0.1cm}
{\bf Proof:}

\vs{0.1cm}
Let $f\in \mathrm{Con}_{R,\mu}(A)$. 
By Remark 4.4 we have to handle the following cases:

\vs{0.2cm}
{\bf Case 1:} $f=\varphi$ where $\varphi\in \mathrm{Sub}_R(A)$.

\vs{0.1cm}
{\bf Proof of Case 1:} This is settled by above.

\vs{0.2cm}
{\bf Case 2:} $f=\log_\mu(\psi)$ where $\psi\in \mathrm{Sub}_R(A)$ is positive.

\vs{0.1cm}
{\bf Proof of Case 2:}
Let $x\in A$.
Then 
\begin{eqnarray*}
f_\Sigma(x)&=&\log_S\big(\psi_\Sigma(\rho(x))\big)\stackrel{\mathrm{Case\;1}}{=}
\log_S\big(\rho(\psi(x))\big)\\
&\stackrel{\Sigma\;\mathrm{companion}}{=}&\rho^Y\big(\log_\mu(\psi(x))\big)=\rho^Y\big(f(x)\big).
\end{eqnarray*}
\hfill$\Box$

\vs{0.5cm}
{\bf 4.6 Corollary}

\vs{0.1cm}
{\it 
The map $\mathfrak{E}_{\mu,\Sigma}^A$ is an embedding.}

\vs{0.1cm}
{\bf Proof:}

\vs{0.1cm}
Let $f,g\in \mathrm{Con}_{R,\mu}(A)$ such that $f_\Sigma=g_\Sigma$. Then $f_\Sigma\circ\rho=g_\Sigma\circ \rho$. By Proposition 4.5 we obtain that $\rho^Y\circ f=\rho^Y\circ g$. Since $\rho^Y$ is an embedding we get that $f=g$.

\hfill$\Box$

\vs{0.5cm}
Given a globally subanalytic subset $B$ of $\IR^n$ we denote by
$\mathfrak{L}_R^B:\mathrm{Sub}_\IR(B)\hookrightarrow \mathrm{Sub}_R(B_R)$
and
$\mathfrak{L}_S^B:\mathrm{Con}_\IR(B)\hookrightarrow \mathrm{Con}_S(B_S)$
the canonical liftings (cf. Sections 1.3 \& 1.4). Note that $B_S=(B_R)_\Sigma$.

\vs{0.5cm}
{\bf 4.7 Theorem}

\vs{0.1cm}
{\it Let $B$ be a globally subanalytic subset of $\IR^n$. The following holds:
	\begin{itemize}
		\item[(1)] 
	There is exactly one map $\mathfrak{L}_{R,\mu}^B:\mathrm{Con}_\IR(B)\to \mathrm{Con}_{R,\mu}(B_R),f\mapsto f_{R,\mu},$ such that the following diagram commutes:
	\begin{center}
		\begin{xy}\;\;\;\;\;\;\;\;\;\;\;\;\;\;\;\;\;\;\;\;\;\;\;\;\;\;\;\;\;\;\;\;\;\;\;\;\;\;\;\;\;\;\;\;\;\;\;\;\;\;\;\;\;\;\;\;\;\;\;\;
			\xymatrix{
				\mathrm{Con}_{R,\mu}(B_R)\ar[rr]^{\mathfrak{E}_{\mu,\Sigma}^{B_R}}& &\mathrm{Con}_S(B_S)\\
				&\mathrm{Con}_\IR(B).\ar[lu]^{\mathfrak{L}_{R,\mu}^B}\ar[ru]_{\mathfrak{L}_S^B}&
			}
		\end{xy}
	\end{center}
\item[(2)] The map $\mathfrak{L}_{R,\mu}^B$ is an embedding that extends $\mathfrak{L}_R^B$. 
\item[(3)] The map $\mathfrak{L}_{R,\mu}^B$ does not depend on the choice of $\Sigma$.
\end{itemize}}
{\bf Proof:}

\vs{0.1cm}
(1): {\bf Uniqueness:}
Let $f\in\mathrm{Con}_\IR(B)$ and let $F,\widetilde{F}\in \mathrm{Con}_{R,\mu}(B_R)$ such that $F_\Sigma=\widetilde{F}_\Sigma=f_S$.
Then by Corollary 4.6 $F=\widetilde{F}$.

\vs{0.2cm}
{\bf Existence:}
Let $f\in \mathrm{Con}_\IR(B)$. Then there are $p,q\in\IN_0$ and for $i\in\{1,\ldots,p\},j\in\{1,\ldots,q\}$ there are $\varphi_i\in\mathrm{Sub}_\IR(B)$ and positive $\psi_{ij}\in\mathrm{Sub}_\IR(B)$ such that 
$$f=\sum_{i=1}^k\varphi_i\Big(\prod_{j=1}^l\log\big(\psi_{ij}\big)\Big).$$
We set
$$f_{R,\mu}:=\sum_{i=1}^k\big(\varphi_i\big)_R\Big(\prod_{j=1}^l\log_\mu\big(\big(\psi_{ij}\big)_R\big)\Big).$$
Then clearly $\big(f_{R,\mu}\big)_\Sigma=f_S$.

\vs{0.2cm}
(2): That $\mathfrak{L}_{R,\mu}$ is an embedding follows from the fact that $\mathfrak{L}^B_S$ is an embedding. That it extends $\mathfrak{L}_R$ follows by the construction in (1).

\vs{0.2cm}
(3): This follows from the proof of (1).
\hfill$\Box$

\vs{0.5cm}
By $\mathrm{res}$ we denote the restriction of functions.

\vs{0.5cm}
{\bf 4.8 Remark}

\vs{0.1cm}
Let $B,C$ be globally subanalytic subsets of $\IR^n$ with $B\subset C$.
Then the following diagram commutes:

\begin{center}
	\begin{xy}\;\;\;\;\;\;\;\;\;\;\;\;\;\;\;\;\;\;\;\;\;\;\;\;\;\;\;\;\;\;\;\;\;\;\;\;\;\;\;\;\;\;\;\;\;\;\;\;\;\;\;\;\;\;\;\;\;\;\;\;
		\xymatrix{
			\mathrm{Con}_{R,\mu}(C_R)\ar[r]^{\mathrm{res}}& \mathrm{Con}_{R,\mu}(B_R)\\
			\mathrm{Con}_\IR(C)\ar[u]^{\mathfrak{L}_{R,\mu}^C}\ar[r]^{\mathrm{res}}& \mathrm{Con}_\IR(B).\ar[u]_{\mathfrak{L}_{R,\mu}^B}
		}
	\end{xy}
\end{center}

We compare the embedding obtained by distinct logarithmic data. 

\vs{0.5cm}
\newpage
{\bf 4.9 Proposition}

\vs{0.1cm}
{\it Let $\mu,\mu'$ be logarithmic data for $R$.
The following are equivalent:
\begin{itemize}
	\item[(i)]  There is a ring homomorphism $\Phi:R[X]\to R[X]$ such that $\Phi\circ \mathfrak{L}_{R,\mu}^B=\mathfrak{L}_{R,\mu'}^B$ for any globally subanalytic subset $B$ of any $\IR^n$.
	\item[(ii)]  There is a ring homomorphism $\Phi:R[X]\to R[X]$ such that $\Phi\circ \mathfrak{L}_{R,\mu}^B=\mathfrak{L}_{R,\mu'}^B$ for any globally subanalytic subset $B$ of some $\IR^n$.
	\item[(iii)] There is a ring homomorphism $\Phi:R[X]\to R[X]$ such that $\Phi\circ \mathfrak{L}_{R,\mu}^{\IR^n}=\mathfrak{L}_{R,\mu'}^{\IR^n}$ for any $n\in\IN$.
	\item[(iv)] There is a ring homomorphism $\Phi:R[X]\to R[X]$ such that $\Phi\circ \mathfrak{L}_{R,\mu}^{\IR^n}=\mathfrak{L}_{R,\mu'}^{\IR^n}$ for some $n\in\IN$.
	\item[(v)] The logarithmic data $\mu$ and $\mu'$ are equivalent. 
\end{itemize}
If one of the above holds then $\Phi$ is the logarithmic connection between $\mu$ and $\mu'$.}

\vs{0.1cm}
{\bf Proof:}

\vs{0.1cm}
The directions (i) $\Rightarrow$ (ii),  (i) $\Rightarrow$ (iii), (ii) $\Rightarrow$ (iv) and (iii) $\Rightarrow$ (iv) are obvious.

\vs{0.1cm}
(iv) $\Rightarrow$ (v):
Let $n$ and $\Phi$ be as in (iv). Let 
$$f:\IR^n\to \IR, x=(x_1,\ldots,x_n)\mapsto\left\{\begin{array}{ccc}
\log(x_1),&&x_1>0,\\
&\mbox{if}&\\
0,&&x_1\leq 0.\\
\end{array}\right.$$ 
Then 
$$f_{R,\nu}:R^n\to R, x=(x_1,\ldots,x_n)\mapsto\left\{\begin{array}{ccc}
\log_\nu(x_1),&&x_1>0,\\
&\mbox{if}&\\
0,&&x_1\leq 0\\
\end{array}\right.$$ 
for $\nu\in\{\mu,\mu'\}$.
Hence $\Phi$ is the logarithmic connection between $\mu$ and $\mu'$.

\vs{0.2cm}
(v) $\Rightarrow$  (i): Let $\Phi:R[X]\stackrel{\cong}{\rightarrow}R[X]$ be the logarithmic connection between $\mu$ and $\mu'$. Let $n\in\IN$ and let $B$ be a globally subanalytic subset of $\IR^n$. 
Let $f\in \mathrm{Con}_\IR(B)$. 
By the proof of Theorem 4.7 and by Proposition 2.10(1) we obtain that 
$\Phi\big(f_{R,\mu}\big)=f_{R,\mu'}$.
\hfill$\Box$

\section{Lebesgue integration}

As an application we show how we can extend the construction in [14, Section 2] to establish on $R$ a Lebesgue measure and integration theory for globally subanalytic sets and functions after fixing a logarithmic datum $\mu$ for $R$.
By $\lambda_n$ we denote the usual Lebesgue measure on $\IR^n$.

\vs{0.5cm}
{\bf 5.1 Construction}

\begin{itemize}
	\item[(I)] 
Let $A$ be a globally subanalytic subset of $R^n$. 
We define its {\bf Lebesgue measure}
$$\lambda_{R,n}^\mu(A)\in R[X]\cup\{\infty\}$$
{\bf with respect to} $\mathbf{\mu}$
as follows. Take an $\ma{L}_\an$-formula $\psi(x,y)$, $x=(x_1,\ldots,x_n),$ $y=(y_1,\ldots,y_q)$, and a point $a\in R^q$ such that $A=\psi(R^n,a)$. 
By Comte, Lion and Rolin [5, 17] the function
$$F:\IR^q\to \IR, y\mapsto \left\{\begin{array}{ccc}
\lambda_n\big(\psi(\IR^n,y)\big),&&\lambda_n\big(\psi(\IR^n,y)\big)<\infty,\\
&\mbox{if}&\\
-1,&&\lambda_n\big(\psi(\IR^n,y)\big)=\infty,\end{array}\right.$$
is constructible.
Set 
$$\lambda_{R,n}^\mu(A):=\left\{\begin{array}{ccc}
F_{R,\mu}(a),&&F_{R,\mu}(a)\geq 0,\\
&\mbox{if}&\\
\infty,&& F_{R,\mu}(a)=-1.\\
\end{array}\right.$$
\item[(II)] 
Let $f:R^n\to R$ be a globally subanalytic function that is non-negative. 
We define its {\bf Lebesgue integral} 
$$\int_{R^n}^\mu f(x)\;dx\in R[X]\cup\{\infty\}$$
{\bf with respect to $\mathbf{\mu}$}
as follows. There is $q\in\IN_0$, a non-negative globally subanalytic function $\psi:\IR^{n+q}\to \IR,(x,y)\mapsto \psi(x,y),$ $x=(x_1,\ldots,x_n),y=(y_1,\ldots,y_q)$, and a point $a\in R^q$ such that $f=\psi_R(-,a)$. 
By Comte, Lion and Rolin [5, 17] the function
$$F:\IR^q\to \IR, y\mapsto \left\{\begin{array}{ccc}
\int_{\IR^n}\psi(x,y)\;dx,&&\int_{\IR^n}\psi(x,y)\;dx<\infty,\\
&\mbox{if}&\\
-1,&&\int_{\IR^n}\psi(x,y)\;dx=\infty,\end{array}\right.$$
is constructible.
Set 
$$\int_{R^n}^\mu f(x)\;dx:=\left\{\begin{array}{ccc}
F_{R,\mu}(a),&&F_{R,\mu}(a)\geq 0,\\
&\mbox{if}&\\
\infty,&& F_{R,\mu}(a)=-1.\\
\end{array}\right.$$
\end{itemize}

{\bf Proof of well-definition:}

\vs{0.1cm}
By Theorem 3.2 there is a $T_{\an,\exp}$-companion $\Sigma=\big(S,\rho,Y\big)$ for $\mu$.
By a routine model theoretic argument, the value $F_S(\rho(a))$ does not depend on the choices of $\psi$ and $a$.
By Theorem 4.7 and Proposition 4.5 we have 
$\rho^Y(F_{R,\mu}(a))=F_S(\rho(a))$.
Since $\rho^Y$ is an embedding the value $F_{R,\mu}(a)$ does not depend on the choices of $\psi$ and $a$.
\hfill$\Box$

\vs{0.5cm}
{\bf 5.2 Remark}

\vs{0.1cm}
As in [14, Section 2] we obtain that the Lebesgue measure with respect to $\mu$ is additive, monotone, translation invariant, fulfills the product formula and extends the naive volume and that the Lebesgue integral (see notations below) is linear and monotone.
The Lebesgue measure and integral with  respect to $\mu$ are intertwined as usually; i.e. the measure is the integral of its characteristic function and the integral of a non-negative function is the measure of its subgraph.

\vs{0.5cm}
Let $A\subset R^n$ be globally subanalytic. We say that $A$ has {\bf finite Lebesgue measure} (with respect to $\mu$) if  $\lambda_{R,n}^\mu(A)<\infty$.\\
Let $f:R^n\to R$ be globally subanalytic. We say that $f$ is {\bf Lebesgue integrable} (with respect to $\mu$) if  $\int_{R^n}^\mu f_+(x)\;dx<\infty$ and $\int_{R^n}^\mu f_-(x)\;dx<\infty$. In that case we set $\int_{R^n}^\mu f(x)\; dx:=\int_{R^n}^\mu f_+(x)\;dx-\int_{R^n} f_-(x)\;dx$.\\
Let $f:R^n\to R$ be globally subanalytic and let $A$ be a globally subanalytic subset of $R^n$. We say that $f$ is {\bf Lebesgue integrable over $A$} (with respect to $\mu$) if $f\mathbbm{1}_A$ is Lebesgue integrable with respect to $\mu$.
We set then $\int_{A}^\mu f(x)\;dx:=\int_{R^n}^\mu f(x)\mathbbm{1}_A(x)\,dx$.
Let $A$ be a globally subanalytic subset of $R^n$ and let $f:A\to R$ be globally subanalytic. Let $\bar{f}$ be the function which extends $f$ by $0$ to $R^n$. We say that $f$ is {\bf Lebesgue integrable} (with respect to $\mu$) if $\bar{f}$ is Lebesgue integrable with respect to $\mu$.
We set then $\int_{A}^\mu f(x)\;dx:=\int_{R^n}^\mu \bar{f}(x)\;dx$.

\vs{0.5cm}
That the notions of finiteness and integrability do indeed not depend on the choice of the logarithmic datum is shown in the next proposition.

\vs{0.5cm}
{\bf 5.3 Proposition}

\vs{0.1cm}
{\it Let $\mu,\mu'$ be logarithmic data for $R$.
	\begin{itemize}
		\item[(1)] Let $A$ be a globally subanalytic subset of $R^n$. Then $A$ has finite Lebesgue measure with respect to $\mu$ if and only if it has finite Lebesgue measure with respect to $\mu'$.
		\item[(2)] Let $f:R^n\to R$ be a globally subanalytic function. Then $f$ is Lebesgue integrable with respect to $\mu$ if and only if it is Lebesgue integrable with respect to $\mu'$. 
	\end{itemize}}
{\bf Proof:}

\vs{0.1cm}
Let $\nu\in \{\mu,\mu'\}$.

\vs{0.2cm}
(1): Let $A$ be as in (1) and let $F$ be as in Construction 5.1. Let $B:=F^{-1}(-1)$. Then $B$ is globally subanalytic by [5].
We obtain by Theorem 4.7(3) and Remark 4.8 that 
$$\lambda_{R,n}^\nu(A)=\infty\Leftrightarrow F_{R,\mu}(a)=-1\Leftrightarrow (\mathbbm{1}_B)_R(a)=1\Leftrightarrow a\in B_R.$$ 
Since the definition of $B_R$ does not depend on the choice of $\nu$ we are done.

\vs{0.2cm}
(2): This is shown in the same way.
\hfill$\Box$

\vs{0.5cm}
Given $n\in\IN$ we denote by $\chi_{R,n}$ the set of globally subanalytic subsets of $R^n$ that have finite Lebesgue measure and by $\ma{L}^1_{R,n}$ the set of integrable globally subanalytic function on $R^n$. We set $\lambda_{R,n}^\mu:\chi_{R,n}\to R[X], A\mapsto \lambda_{R,n}^\mu (A),$ and $\mathrm{Int}_{R,n}^\mu:\ma{L}_{R,n}^\mu(A), f\mapsto \int_{R^n}^\mu f(x)\;dx$.

\vs{0.5cm}
{\bf 5.4 Proposition}

\vs{0.1cm}
{\it Let $\mu,\mu'$ be logarithmic data for $R$.
	The following are equivalent:
	\begin{itemize}
		\item[(i)]  There is a ring homomorphism $\Phi:R[X]\to R[X]$ such that $\Phi\circ \lambda_{R,\mu}^n=\lambda_{R,\mu'}^n$ for any $n\in\IN$.
		\item[(ii)]  There is a ring homomorphism $\Phi:R[X]\to R[X]$ such that $\Phi\circ \lambda_{R,\mu}^n=\lambda_{R,\mu'}^n$ for some $n\in\IN$ with $n\geq 2$.
		\item[(iii)] The logarithmic data $\mu$ and $\mu'$ are equivalent. 
	\end{itemize}
	If one of the above holds then $\Phi$ is the logarithmic connection between $\mu$ and $\mu'$.}

\vs{0.1cm}
{\bf Proof:}

\vs{0.1cm}
The direction (i) $\Rightarrow$ (ii) is obvious.

\vs{0.2cm}
(ii) $\Rightarrow$ (iii): Let $n\geq 2$ and $\Phi$ be as in (ii).
Let $a\in R$ with $a\geq 1$ and set
$$A_a:=\Big\{(x_1,\ldots,x_n)\in R^n\;\Big\vert\; 1\leq x_1\leq a, 0\leq x_2\leq 1/x_1,0\leq x_j\leq 1\mbox{ for }j>2\Big\}.$$ 
Then $A_a$ is a semialgebraic and hence globally subanalytic subset of $R^n$.
We have that $\lambda_{R,n}^\mu(A_a)=\log_\mu(a)$ and $\lambda_{R,n}^{\mu'}(A_a)=\log_{\mu'}(a)$.
Therefore $\Phi\circ \log_\mu|_{[1,\infty[}=\Phi\circ\log_{\mu'}|_{[1,\infty[}$. By the logarithmic rule (see Theorem 2.6) we get that 
$\Phi\circ \log_\mu=\log_{\mu'}$. Hence $\Phi$ is the logarithmic connection between $\mu$ and $\mu'$.

(iii) $\Rightarrow$ (i): Let $\Phi:R[X]\stackrel{\cong}{\rightarrow}R[X]$ be the logarithmic connection between $\mu$ and $\mu'$.
Let $n\in\IN$ and let $A\in \chi_{R,n}$. 
Let $F:\IR^q\to \IR$ and $a\in R^q$ be as in Construction 5.1.
Since $A$ has finite Lebesgue measure we have that
$\lambda_{R,n}^\nu(A)=F_{R,\nu}(a)$ for $\nu\in \{\mu,\mu'\}$. We are done since by Proposition 4.9 (v)$\Rightarrow$(iii) we have that $\Phi\circ F_{R,\mu}=F_{R,\mu'}$.
\hfill$\Box$

\vs{0.5cm}
Note that $n=1$ is not sufficient in Proposition 5.4(ii) since the measure of an interval is its length (see Remark 5.2; see also [10, Theorem 2.3]).
So we do only get values in $R$, the logarithm does not show up.

\vs{0.5cm}
\newpage
{\bf 5.5 Proposition}

\vs{0.1cm}
{\it Let $\mu,\mu'$ be logarithmic data for $R$.
	The following are equivalent:
	\begin{itemize}
		\item[(i)]  There is a ring homomorphism $\Phi:R[X]\to R[X]$ such that $\Phi\circ \mathrm{Int}_{R,n}^\mu=\mathrm{Int}_{R,n}^{\mu'}$ for any $n\in\IN$.
		\item[(ii)]  There is a ring homomorphism $\Phi:R[X]\to R[X]$ such that $\Phi\circ \mathrm{Int}_{R,n}^\mu=\mathrm{Int}_{R,n}^{\mu'}$ for some $n\in\IN$.
		\item[(iii)] The logarithmic data $\mu$ and $\mu'$ are equivalent. 
	\end{itemize}
	If one of the above holds then $\Phi$ is the logarithmic connection between $\mu$ and $\mu'$.}

\vs{0.1cm}
{\bf Proof:}

\vs{0.1cm}
This follows from Proposition 5.4 and Remark 5.2 on the interplay between the measure and the integral.

\vs{0.5cm}
In view of Propositions 5.4 \& 5.5 and Theorem 2.16 we say in the case that the rational rank coincides with the archimedean rank of the value group of $R$ that the Lebesgue measure and integration theory on $R$ is {\bf uniquely determined up to unique isomorphisms}.

\vs{0.5cm}
{\bf 5.6 Proposition}

\vs{0.1cm}
{\it Let $f:R^m\times R^n\to R,(x,y)\mapsto f(x,y),$ be globally subanalytic.
	\begin{itemize}
		\item[(i)]
	The set
	$$\mathrm{Fin}(f):=\Big\{x\in R^m\;\Big\vert\;f(x,-)\mbox{ is Lebesgue integrable}\Big\}$$
	is globally subanalytic.
	\item[(ii)]
	The function 
	$$\mathrm{Fin}(f)\to R,x\mapsto \int_{R^n}^\mu f(x,y)\;dy,$$
	is $\mu$-constructible.
\end{itemize}} 
{\bf Proof:}

\vs{0.1cm}
This follows by doing Construction 5.1 with parameters and by using the results of [5].
\hfill$\Box$

\vs{0.5cm}
As in [14, Section 5] one can now establish the main theorems of integration in the globally subanalytic context: the transformation formula (see [14, Theorem 5.1]), Lebesgue's theorem on dominated convergence (see [14, Theorem 5.2]) and the fundamental theorem of calculus (see [14, Theorem 5.7]). To obtain Fubini's theorem we have to extend the integration to constructible functions as follows.

\vs{0.5cm}
{\bf 5.7 Construction}

\vs{0.1cm}
Let $n\in\IN$ and  let $f:R^n\to R[X]$ be $\mu$-constructible. 
We define that $f$ is {\bf Lebesgue integrable with respect to $\mu$} and in that case its {\bf Lebesgue integral}
$$\int_{R^n}^\mu f(x)\;dx\in R[X]$$
{\bf with respect to $\mathbf{\mu}$} as follows. There is $q\in \IN_0$, a constructible function $\psi:\IR^{n+q}\to \IR, (x,y)\mapsto \psi(x,y)$, $x=(x_1,\ldots,x_n),y=(y_1,\ldots,y_q)$, and a point $a\in R^q$ such that $f=\psi_{R,\mu}(-,a)$. 
By Cluckers and Dan Miller [1, 2, 3] there are constructible functions $g,h:\IR^q\to \IR$ such that $\psi(-,y)$ is Lebesgue integrable
if and only if $h(y)=0$ and
in that case that $\int_{\IR^n}\psi(x,y)\,dx=g(y)$. 
Then $f$ is Lebesgue integrable with respect to $\mu$ if $h_{R,\mu}(a)= 0$ and we set in this case 
$$\int_{R^n}^\mu f(x)\;dx:=g_{R,\mu}(a).$$

\vs{0.1cm}
{\bf Proof of well-definition:}

\vs{0.1cm}
By Theorem 3.2 there is a $T_{\an,\exp}$-companion $\Sigma=\big(S,\rho,Y\big)$ for $\mu$.
By a routine model theoretic argument, the values $g_S(\rho(a))$ and $h_S(\rho(a))$ do not depend on the choices of $\psi$ and $a$.
By Theorem 4.7 and Proposition 4.5 we have 
$\rho^Y(g_{R,\mu}(a))=g_S(\rho(a))$ and $\rho^Y(h_{R,\mu}(a))=h_S(\rho(a))$.
Since $\rho^Y$ is an embedding the values $g_{R,\mu}(a)$ and $h_{R,\mu}(a)$ do not depend on the choices of $\psi$ and $a$.
\hfill$\Box$

\vs{0.5cm}
{\bf 5.8 Proposition}

\vs{0.1cm}
{\it Let $f:R^m\times R^n,(x,y)\mapsto f(x,y),$ be $\mu$-constructible.
	Then there are $\mu$-constructible functions $g,h:R^m\to R$ such that
	$f(x,-)$ is Lebesgue integrable with respect to $\mu$ if and only if $h(x)= 0$ and that
	$$\int_{R^n}^\mu f(x,y)\;dy=g(x)$$ for such an $x$.}

\vs{0.1cm}
{\bf Proof:}

\vs{0.1cm}
This follows by doing Construction 5.7 with parameters.
\hfill$\Box$

\vs{0.5cm}
In this setting, Fubini's theorem (cf. [14, Theorem 5.12]) can be established.

\vs{0.5cm}
As an application of the Lebesgue theory we obtain as in [14, Section 6] the following approximation result.

\vs{0.5cm}
\newpage
{\bf 5.9 Theorem}

\vs{0.1cm}
{\it Let $a,b\in R$ with $a<b$ and let $f:[a,b]\to R$ be globally subanalytic and continuous. Let $\varepsilon\in R_{>0}$. Then there is some open interval $I$ in $R$ containing $[a,b]$ and some globally subanalytic function $u:I\to R$ that is $C^\infty$ such that $|f(x)-u(x)|<\varepsilon$ for all $x\in [a,b]$.}

\vs{0.5cm}
{\bf 5.10 Final remark}

\begin{itemize}
	\item[(1)] In [14] we have worked in the setting of real closed fields containing the reals with archimedean value group and have obtained a Lebesgue measure and integration theory for the semialgebraic sets and functions.
	Given a real closed field $R$ containing the reals such that its value group $\Gamma$ has finite archimedean rank $\ell$ we can define with the results of the present paper a Lebesgue measure and integration theory for the semialgebraic category by an 
	embedding $R\hookrightarrow \IR((t^\Gamma))$ (cf. Fact 1.2). But we do not obtain all results of [14, Section 3] in the semialgebraic setting since [14, Proposition 1.6] cannot be generalized beyond the case of an archimedean value group.
\item[(2)] Given a model $R$ of the theory $T_\an$ with archimedean value group $\Gamma$ we have constructed in [14, Section 3.4] a Lebesgue measure and integration theory (called there analytic) for the globally subanalytic category.
In [14] the construction relies on a so-called Lebesgue datum $\alpha:=\big(s,\sigma,\tau\big)$ where $\big(s,\tau\big)$ is a logarithmic datum in the sense of the present paper and $\sigma:R\hookrightarrow \IR((t^\Gamma))$ is an embedding compatible with the section $s$. In [14] we have worked with the special $\IR_{\an,\exp}$-companion $\big(\Theta_\alpha,\IR((t))^\mathrm{LE},\log(t^{-1})\big)$ where $\IR((t))^\mathrm{LE}$ is the field of logarithmic-exponential series and the $\ma{L}_\an$-embedding $\Theta_\alpha:R\hookrightarrow \IR((t))^\mathrm{LE}$ is given by
$$R\stackrel{\sigma}{\hookrightarrow} \IR((t^\Gamma))\hookrightarrow\IR((t^\IR))\stackrel{\mathrm{incl}}{\hookrightarrow}\IR((t))^\mathrm{LE}$$
where the map in the middle is defined as $\sum a_\gamma t^\gamma\mapsto \sum a_\gamma t^{\tau(\gamma)}$.
So even in the case $\ell=1$ the present setting is more general than the one in [14].
Moreover, in [14] the measure respectively integral takes in general values in $\IR((t^\Gamma))[X]$. In this paper the construction has been set up such that the values are in $R[X]$.
\end{itemize}

\vs{1cm}
\newpage
\noi {\footnotesize \centerline{\bf References}
	\begin{itemize}
		\item[(1)] 
		R. Cluckers and D. Miller:
		Stability under integration of sums of products of real globally subanalytic functions and their logarithms.
		{\it Duke Math. J.} {\bf 156} (2011), no. 2, 311-348.
		\item[(2)] 
		R. Cluckers and D. Miller:
		Loci of integrability, zero loci, and stability under integration for constructible functions on Euclidean space with Lebesgue measure.
		{\it Int. Math. Res. Not.} (2012), no. 14, 3182-3191.
		\item[(3)]
		R. Cluckers and D. Miller:
		Lebesgue classes and preparation of real constructible functions.
		{\it J. Funct. Anal.} {\bf 264} (2013), no. 7, 1599-1642.
		\item[(4)] R. Cluckers, J. Nicaise and J. Sebag:
	    Motivic integration and its interactions with model theory and Non-Archimedean geometry, Volume I.
        {\it London Math. Soc. Lecture Note Ser.} {\bf 383}, Cambridge Univ. Press, 2011.
		\item[(5)] 
		G. Comte, J.-M. Lion and J.-P. Rolin:
		Nature log-analytique du volume des sous-analytiques.
		{\it Illinois J. Math.} {\bf 44} (2000), no. 4, 884-888.
		\item[(6)] 
		L. van den Dries, A. Macintyre and D. Marker:
		The elementary theory of restricted analytic fields with exponentiation.
		{\it Ann. Math.} {\bf 140} (1994), 183-205.
		\item[(7)] 
		L. van den Dries, A. Macintyre and D. Marker: Logarithmic-exponential series.
		{\it Ann. Pure Appl. Logic} {\bf 111} (2001), no. 1-2, 61-113.
		\item[(8)] 
		A, Ducros, E. Hrushovski and F. Loeser:
		Non-archimedean integrals as limits of complex integrals.
		{\it arXiv:1912.09162v2} (2020), 55 pp.
		\item[(9)] 
		W. Hodges:
		A shorter model theory.
		Cambridge University Press, 1997.
		\item[(10)] 
		T. Kaiser:
		First order tameness of measures.
		{\it Ann. Pure Appl. Logic} {\bf 163} (2012), no. 12, 1903-1927.
		\item[(11)] 
		T. Kaiser:
		Integration of semialgebraic functions and integrated Nash functions.
		{\it Math. Z.} {\bf 275} (2013), no. 1-2, 349-366.
		\item[(12)]
		T. Kaiser: Global complexification of real analytic globally subanalytic functions.
		{\it Israel J. of Math.} {\bf 213} (2016), no. 1, 139-174.
		\item[(13)]
		T. Kaiser: R-analytic functions.
		{\it Arch. Math. Log.} {\bf 55} (2016), no. 5-6, 605-623.
		\item[(14)] 
		 T. Kaiser: Lebesgue measure and integration theory on non-archimedean real closed fields with archimedean value group.
		{\it Proc. Lond. Math. Soc.} {\bf 116} (2018), no. 2, 209-247. 
	    \item[(15)] 
		F.-V. Kuhlmann, S. Kuhlmann and S. Shelah:
		Exponentiation in power series fields.
		{\it Proc. Amer. Math. Soc.} {\bf 125} (1997), no. 11, 3177-3183.
		\item[(16)] S. Kuhlmann:
	Ordered exponential fields.
{\it Fields Institute Monographs}, American Mathematical Society, 2000.
		\item[(17)] 
		J.-M. Lion and J.-P. Rolin:
		Int\'{e}gration des fonctions sous-analytiques et volumes des sous-ensembles sous-analytiques.
		{\it Ann. Inst. Fourier (Grenoble)} {\bf 48} (1998), no. 3, 755-767.
		\item[(19)] 
		S. Prie{\ss}-Crampe:
		Angeordnete Strukturen: Gruppen, K\"orper, projektive Ebenen.
		Springer, 1983.
\end{itemize}}

\vs{0.5cm}
Tobias Kaiser\\
University of Passau\\
Faculty of Computer Science and Mathematics\\
tobias.kaiser@uni-passau.de\\
D-94030 Germany

\end{document}